\documentclass[11pt,reqno]{amsart}

\newcommand{\mysection}[1]{\section{#1}
      \setcounter{equation}{0}}

\newcommand{\nlimsup}{\operatornamewithlimits{\overline{lim}}}

\newtheorem{theorem}{Theorem}[section]
\newtheorem{lemma}[theorem]{Lemma}

\newtheorem{corollary}[theorem]{Corollary}

\theoremstyle{definition}
\newtheorem{assumption}{Assumption}[section]

\theoremstyle{remark}
\newtheorem{remark}{Remark}[section]

 \makeatletter
 \def\dashint{%
 \operatorname%
 {\,\,\text{\bf--}\kern-.98em\DOTSI\intop\ilimits@\!\!}}
 \makeatother

\newcommand{\WO}{\overset{\scriptscriptstyle0}%
{\mathcal W}\,\!}

\newcommand\bR{\mathbb{R}}

\newcommand\cL{\mathcal{L}}

\newcommand\cW{\mathcal{W}}

\begin{document}

\title[Equations
with growing coefficients]{
On   linear
elliptic and parabolic equations
with growing drift in Sobolev spaces without weights}

\author[N.  Krylov]{N.V. Krylov}%
\thanks{The work   was partially supported
  by NSF grant DMS-0653121}
\address{127 Vincent Hall, University of Minnesota,
Minneapolis,
       MN, 55455, USA}
\email{krylov@math.umn.edu}


\begin{abstract}
We consider uniformly elliptic and parabolic second-order equations with
bounded zeroth-order and bounded VMO leading coefficients and
possibly growing first-order coefficients. We look for solutions
  which are summable to the $p$-th power with respect to the usual
Lebesgue measure along with their first and second-order derivatives
with respect to the spatial variable.
\end{abstract}

\maketitle

\mysection{Introduction}

In this paper we concentrate on problems in the whole space
for uniformly elliptic and parabolic second-order equations with
bounded leading and zeroth-order coefficients and
possibly growing first-order coefficients. We look for solutions
  which are summable to the $p$-th power with respect to the usual
Lebesgue measure along with their first- and second-order derivatives
with respect to the spatial variables.

There exists a quite extensive literature related to equations
with growing coefficients in Sobolev-Hilbert spaces with weights.
Since here no weights are used we only
  refer the reader to \cite{BL}, \cite{CV}, \cite{GL}, \cite{Gy93},
\cite{GK} where one can find further references as well.

 It is generally believed that introducing weights is the 
most natural setting for equations with growing coefficients.
 The present paper seems to be the first
one treating the unique solvability of
these equations in Sobolev spaces $W^{2}_{p}$ for $p\in(1,\infty)$
 without weights and
without imposing any special conditions on the relations between
the coefficients or on their smoothness. In the elliptic case,
in rough terms,
 it is sufficient for us that the drift term $b^{i}D_{i}u$
be, say, such that
\begin{equation}
                                                         \label{2.21.1}
\lim_{\alpha\downarrow0}\sup_{x,y:|x-y|\leq\alpha}|x-y|^{\varepsilon}
|b(x)-b(y)|=0,
\end{equation}
where (possibly negative) $\varepsilon<(d-1)/(d\vee p)$. This condition
has nothing to do with any continuity property of $b$
since $\varepsilon$ is allowed to be positive.

It is worth noting that many issues for 
divergence-type equations with time independent 
growing coefficients in $L_{p}$ spaces without weights
were treated previously in the literature. 
This was done mostly  by using the semigroup approach.
We briefly mention only a few recent papers
sending the reader to them for additional references.

In \cite{LV} a strongly continuous
in $L_{p}$ semigroup is constructed corresponding
to elliptic  operators with measurable
leading coefficients and Lipschitz
continuous drift coefficients. This did not lead to the solvability
of elliptic equations in $W^{1}_{p}$ for $p>2$ because of low regularity
of the leading coefficients.
In \cite{MP} it is assumed that
 if, for $| x|\to\infty$, the drift coefficient   grows,
then  the zeroth-order coefficient  should grow, basically,
as the square of the drift. There  is also a condition on the divergence
of the drift coefficient.
In \cite{PR} there is no zeroth-order term
and the semigroup is constructed under some assumptions
one of which translates into  the monotonicity of 
$\pm b(x)-Kx$, for a constant $K$, if the leading term is the Laplacian.
In \cite{CF}  the drift coefficient
is assumed to be globally Lipschitz
continuous if the zeroth-order coefficient is constant.

Some conclusions in the above cited papers are quite similar to ours
but the corresponding assumptions are not as general
in what concerns the regularity of the coefficients.
However, these papers contain a lot of additional
important information not touched upon in the present paper
(in particular, it is shown in \cite{LV} that the corresponding semigroup
is not analytic).

The technique, we apply, originated from \cite{KP} and uses special
 cut-off functions whose support evolves in time
in a manner adapted to the drift.
Another less important feature is that the leading coefficients
of the equations are assumed to be only measurable in time and VMO
in $x$. In fact, the reader will see from our proofs that
nothing special is needed from the leading terms and 
one can add the drift term satisfying our conditions
to any equation for which the Sobolev space theory is available.
In particular,
this can be done for divergence form equations
with measurable coefficients if $p=2$. However, for the sake
of brevity and clarity we concentrate only on nondivergence
type equations.  The main emphasis here is that we allow $b(t,x)$ to grow
as $|x|\to\infty$ and still measure the size of the second-order
derivatives with respect to Lebesgue measure
thus avoiding using weights.

Let $\bR^{d}$ be a Euclidean space of points $x=(x^{1},...,x^{d})$.
 We consider the following second-order
operator $L$:
\begin{equation}                                    \label{lu}
 L u (t,x) = a^{ij}(t,x)D_{ij}u (t,x) + b^{i}(t,x)
D_{i}u (t,x) -c(t,x)u(t,x),
\end{equation}
   acting on functions defined
on $\bR^{d+1}_{T}$, which is
 $[T, \infty) \times \bR^d$ if $T \in (-\infty,  \infty)$ and
   on $\bR^{d+1}$ if $T=-\infty$
(the summation convention is enforced throughout the article).
Here
$$
D_{i}=\frac{\partial}{\partial x^{i}},\quad
D_{ij}=D_{i}D_{j}.
$$
   We are dealing with the parabolic equation
\begin{equation}
                                                \label{2.6.4}
\partial_{t}u (t,x) + Lu(t,x) = f(t,x),\quad (t,x) \in ( T,   \infty)
\times \bR^{d},
\end{equation}
where $\partial_{t}=\partial/\partial t$,
and, in case the coefficients are independent of $t$, with
the elliptic equation
\begin{equation}
                                                \label{2.6.5}
  Lu( x) = f( x),\quad  x \in \bR^{d}.
\end{equation}
The solutions of \eqref{2.6.5} are sought in $W^{2}_{p}(\bR^{d})$,
usual Sobolev space, and the space of solutions
of \eqref{2.6.4} will be $\cW^{2}_{p}(T)$ which 
is defined as follows.

We write $u\in\cW^{2}_{p}(T)$ if $u=u(t,x)$ is a (measurable)
function defined on $\bR^{d+1}_{T}$ such that
\begin{equation}
                                                \label{2.6.6}
\|u\|_{\cL_{p}(\bR^{d+1}_{T})}+
\|Du\|_{\cL_{p}(\bR^{d+1}_{T})}+\|D^{2}u\|_{\cL_{p}
(\bR^{d+1}_{T})}<\infty
\end{equation}
and $\partial_{t}u:=\partial u/\partial t$ is
 locally summable on $\bR^{d+1}_{T}$. 
Of course, $Du$ and $D^{2}u$ are the gradient 
and the Hessian matrix of $u$, respectively.
Observe that we do not include
$\partial_{t}u $ into the left-hand side of \eqref{2.6.6}
because we believe that, generally, in our situation
$\partial_{t}u\not\in \cL_{p}(\bR^{d+1}_{T})$
(see Remark \ref{remark 2.11.1}).

Our main results are presented in Sections
\ref{section 2.15.1} (elliptic case)
and \ref{section 2.15.2} (parabolic case).
Theorem \ref{theorem 2.7.1} saying that under appropriate conditions
the elliptic equation $Lu-\lambda u=f$ is uniquely solvable
in $W^{2}_{p}(\bR^{d})$ if $\lambda$ is large enough
is proved in Section \ref{section 2.15.2}.
Interestingly enough, even if $b$ is constant we do not know
any other proof of Theorem \ref{theorem 2.7.1} not using
the parabolic theory.

We prove Theorems \ref{theorem 2.4.2}
and \ref{theorem 2.9.1}  in Section \ref{section 2.15.4} 
and  \ref{section 2.15.5},  respectively, after we prepare
necessary  tools in Section \ref{section 2.15.3}.
In Section \ref{section 2.15.6} we give an example showing
that for elliptic equations one cannot take $\lambda_{0}>0$
arbitrary small in contrast with the case of bounded coefficients
as described in Section 11.6 of \cite{book}. This fact is known
from \cite{Me},
where the spectrum of the Ornstein-Uhlenbeck operator is found
in the multidimensional case in $\cL_{p}$ spaces and it is shown that
the spectrum depends on $p$.

As usual when we speak of
 ``a constant" we always mean ``a finite constant".

The author is sincerely grateful to A. Lunardi for the fruitful
discussion of the results.

\mysection{Main result for elliptic case}
                                               \label{section 2.15.1}

For $p\in(1,\infty)$, $p\ne d$, define
$$
q=d\vee p,
$$
and if $p=d$ let $q$ be a fixed number such that $q>d$.

\begin{assumption} 
                                       \label{assumption 2.7.2}
(i) The functions $a^{ij},b^{i},c$ are measurable,
$a^{ij}=a^{ji}$, $c\geq 0$.

(ii) There exist   constants $K,\delta>0$ such that
for all values of arguments and $\xi\in\bR^{d}$
$$
\delta|\xi|^{2}\leq a^{ij}\xi^{i}\xi^{j}\leq K|\xi|^{2},
\quad c\leq K.
$$

(iii) The function   $|b |^{q}$
 is locally integrable on $\bR^{d}$.
\end{assumption}
 
The following assumptions contain    parameters $\gamma_{a},
\gamma_{b}\in(0,1]$
whose value will be specified later.
For $\alpha>0$ we denote $B_{\alpha}=\{x\in\bR^{d}:|x|<\alpha\}$.

 \begin{assumption}[$\gamma_{b}$]
                                       \label{assumption 2.3.1}
There exists an $\alpha\in(0,1]$ such that on $\bR^{d }$  
\begin{equation}
                                                        \label{2.15.1}
\alpha^{-d}\int_{B_{\alpha}}\int_{B_{\alpha}}|b(  x+y)
-b(  x+z)|^{q}\,dydz  \leq\gamma_{b}.
\end{equation}
 
\end{assumption}
It is easy to check that
Assumption \ref{assumption 2.3.1} is satisfied
with any $\gamma_{b}>0$ if
\eqref{2.21.1} holds. For instance, we allow $b$
such that $|b(x)-b(y)|\leq K$ if $|x-y|\leq1$.
 We see that
$|b(x)|$ can grow to infinity as $|x|\to\infty$.

 \begin{assumption}[$\gamma_{a}$]
                                       \label{assumption 2.4.01}
There exists an $\varepsilon_{0}>0$ such that
for any $\varepsilon\in(0,\varepsilon_{0}]$, $x\in\bR^{d}$,
 and $i,j=1,...,d$  we have 
\begin{equation}
                                                        \label{2.15.2}
\varepsilon^{-2d } 
 \int_{B_{\varepsilon}}\int_{B_{\varepsilon}}
|a^{ij}( x+y)-a^{ij}( x+z)|\,dydz \leq\gamma_{a}.
\end{equation}

\end{assumption}

Obviously, the left-hand side of \eqref{2.15.2}
is less than
$$
N(d)\sup_{|x-y|\leq2\varepsilon}|a^{ij}( x )-a^{ij}(y)|,
$$
which implies that Assumption \ref{assumption 2.4.01}
 is satisfied with any $\gamma_{a}>0$ if, for instance,
$a$ is a uniformly continuous function.
Recall that if Assumption \ref{assumption 2.4.01}
 is satisfied with any $\gamma_{a}>0$, then one says that
$a$ is in VMO.

Here is one of the main results of the paper.
\begin{theorem}
                                             \label{theorem 2.7.1}
There exist   constants  
$$
\gamma_{a}=\gamma_{ a}(d,\delta,K,p)>0,\quad
\gamma_{b}=\gamma_{ b}(d,\delta,K,p,\varepsilon_{0})>0,
$$
$$
 N=N(d,\delta,K,p,\varepsilon_{0}),
\quad \lambda_{0}=\lambda_{0}(d,\delta,K,p,\varepsilon_{0}
,\alpha )\geq0
$$ 
such that, if the above assumptions are satisfied, then
 for any   $u \in W^{ 2}_{p}(\bR^{d})$ 
and $\lambda\geq\lambda_{0}$ 
  we have
\begin{equation}
                                             \label{2.4.01}
\lambda\|u\|_{\cL_{p}(\bR^{d})}+\|D^{2}u\|_{\cL_{p}(\bR^{d})}
\leq N\|L u -\lambda u\|_{ \cL_{p}(\bR^{d})}.
\end{equation}
Furthermore, for any $f\in \cL_{p}(\bR^{d})$
and $\lambda\geq\lambda_{0}$ there is a unique
$u \in W^{ 2}_{p}(\bR^{d})$ such that $Lu-\lambda u=f$.
\end{theorem}

We prove this theorem in Section \ref{section 2.15.2}.
One of surprising features of \eqref{2.4.01} is that $N$
is independent of $b$ if $b$ is constant. Another one is that
the set $(L-\lambda)W^{ 2}_{p}(\bR^{d})$ may not coincide
with $\cL_{p}(\bR^{d})$ if $|b|$ grows and yet it always
contains $\cL_{p}(\bR^{d})$.  Some consequences of this
peculiarity are discussed in \cite{KP}.

\mysection{Main results for parabolic case}
                                               \label{section 2.15.2}
\begin{assumption} 
                                       \label{assumption 2.7.1}
(i) Assumptions \ref{assumption 2.7.2} (i) and (ii)
are satisfied.

(ii)  
 For any $x\in\bR^{d}$ and $\alpha\in(0,1]$ the function 
$$
 \int_{B_{\alpha}}|b(t,x+y)|\,dy
$$
is locally integrable to the power $p/(p-1)$ with respect to
$t$.
\end{assumption} 

Notice that a simple covering argument shows that
for any $\alpha\in(0,\infty)$ the function 
$$
\sup_{|x|\leq\alpha}\int_{B_{\alpha}}|b(t,x+y)|\,dy
$$
is also locally integrable to the power $p/(p-1)$ with respect to
$t$.

 \begin{assumption}[$\gamma_{b}$]
                                       \label{assumption 2.8.1} 
There exists an $\alpha\in(0,1]$ such that on $\bR^{d+1}$ (a.e.)
$$
\alpha^{-d}\int_{B_{\alpha}}\int_{B_{\alpha}}|b(t, x+y)
-b(t, x+z)|^{q}\,dydz  \leq\gamma_{b}.
$$

\end{assumption}

 \begin{assumption}[$\gamma_{a}$]
                                       \label{assumption 2.4.1} 
There  exists an $\varepsilon_{0}>0$ such that
for any $\varepsilon\in(0,\varepsilon_{0}]$,  
 $s\in \bR$, and $i,j=1,...,d$, we have
\begin{equation}
                                             \label{2.9.5}
\varepsilon^{-2d-2}\int_{s}^{s+\varepsilon^{2}}\bigg(
 \sup_{x\in \bR^{d}}\int_{B_{\varepsilon}}\int_{B_{\varepsilon}}
|a^{ij}(t,x+y)-a^{ij}(t,x+z)|\,dydz\bigg)
 \,dt\leq\gamma_{a}.
\end{equation}

\end{assumption}

The following is a parabolic analog of 
the estimate in Theorem  \ref{theorem 2.7.1}.
\begin{theorem}
                                             \label{theorem 2.4.2}
There exist   constants  
$$
\gamma_{a}=\gamma(d,\delta,K,p)>0,\quad
\gamma_{b}=\gamma(d,\delta,K,p,\varepsilon_{0})>0,
$$
$$
 N=N(d,\delta,K,p,\varepsilon_{0} ),
\quad \lambda_{0}=\lambda_{0}(d,\delta,K,p,\varepsilon_{0},
\alpha )\geq0
$$  
such that, if the above assumptions are satisfied, then
 for any $T\in[-\infty,\infty)$,  $u \in 
\cW^{ 2}_{p}(T)$,
and $\lambda\geq\lambda_{0}$ 
  we have
\begin{equation}
                                             \label{2.4.1}
\lambda\|u\|_{\cL_{p}(\bR^{d+1}_{T})}+
\|D^{2}u\|_{\cL_{p}(\bR^{d+1}_{T})}
\leq N\|L u+\partial_{t}u -\lambda u\|_{\cL_{p}(\bR^{d+1}_{T})}.
\end{equation}

\end{theorem}
Observe that, if the right-hand side of \eqref{2.4.1}
is finite, then $\partial_{t}u+b^{i}D_{i}u$
is  in $\cL_{p}(\bR^{d+1}_{T})$ and, since $\partial_{t}u$
is locally summable, the same is true for $b^{i}D_{i}u$.
Therefore, not surprisingly,
to prove the existence of solutions
of parabolic equations we impose one more assumption
on $b$,
which would guarantee that $b^{i}D_{i}u$
is locally summable if $u\in\cW^{2}_{p}(T)$.
For $p<d$ set
$$
q_{1}=\frac{pd}{(p-1)d+p},\quad r_{1}=\frac{(p-1)d+p}
{(p-1)d},
$$
and for $p\geq d$ let $q_{1}\in(1,p)$ be any fixed number
and 
$$
r_{1}=\frac{p}{(p-1)q_{1}}.
$$
Observe that $1<q_{1}<q$.
\begin{assumption}
                                        \label{assumption 2.9.1}
For any $s,t\in\bR$, such that $s<t$, and $R\in(0,\infty)$ we have
\begin{equation}
                                             \label{2.10.1}
\int_{s}^{t}\big(\int_{B_{R}}|b(\tau,x)|^{q_{1}}
\,dx\big)^{r_{1}}\,d\tau<\infty.
\end{equation}
\end{assumption}
Notice that this assumption
coinsides with Assumption \ref{assumption 2.7.1} (ii)
if $b$ is independent of $x$.

\begin{theorem}
                                             \label{theorem 2.9.1}
Take the  constants  
$ 
\gamma_{a} ,
\gamma_{b}$, and
 $
 \lambda_{0} 
$  
from Theorem \ref{theorem 2.4.2} and suppose that
  Assumptions \ref{assumption 2.7.1}-\ref{assumption 2.9.1}
are satisfied. Then
 for any $\lambda\geq\lambda_{0}$,
$f\in \cL_{p}(\bR^{d+1}_{T})$,
and $T\in[-\infty,\infty)$, there is a unique
$u \in \cW^{ 2}_{p}(T)$ such that $
\partial_{t}u+Lu-\lambda u=f$ in $\bR^{d+1}_{T}$.
 
\end{theorem}
We also have a result for the Cauchy problem.
Fix $T,S\in\bR$ such that $T<S$ and write
$u\in\WO^{2}_{p}(T,S)$ if $u\in\cW^{2}_{p}(T)$
and $u(t,x)=0$ for $t\geq S$.

\begin{theorem}
                                       \label{theorem 2.7.3}
Take the  constants  
$ \gamma_{a}$, and
$\gamma_{b}$ 
from Theorem \ref{theorem 2.4.2} and suppose that
  Assumptions \ref{assumption 2.7.1}-\ref{assumption 2.9.1}
are satisfied. Then for any 
$$
f\in \cL_{p}((T,S)\times\bR^{d } ),\quad
v\in W^{1,2}_{p}((T,\infty)\times\bR^{d } )
$$
 there exists
a unique $u\in\cW^{2}_{p}(T)$ such that $\partial_{t}u
+Lu=f$ in $(T,S)\times\bR^{d } $ and $u-v\in 
\WO^{2}_{p}(T,S)$

\end{theorem}

Proof. It suffices to prove the theorem for the equation
$\partial_{t}u
+Lu-\lambda u=f$ with $\lambda$ as large as we like.
We take it so large that we can apply Theorem \ref{theorem 2.9.1}.
Next, we change the coefficients of $L$ for $t\geq S$
 if needed in such a way
that $L=\Delta$ for $t\geq S$. Finally, we change $f$
 for $t\geq S$
 if necessary and set it to be $(\partial_{t}+\Delta
-\lambda)(\zeta v)$ for $t\geq S$, where $\zeta(t)$
is any $C^{\infty}_{0}(\bR)$ function such that $\zeta(t)=1$
for $t\in(T,S)$. With this new objects 
according to Theorem \ref{theorem 2.9.1}
we can find 
a  $\tilde{u} \in \cW^{ 2}_{p}(T)$ such that $
\partial_{t}\tilde{u}+L\tilde{u}-\lambda 
 \tilde{u}=f$ in $\bR^{d+1}_{T}$.
After applying Theorem \ref{theorem 2.4.2} with $S$ and $\tilde{u}
-\zeta v$
in place of $T$ and $u$, respectively, we see that
$\tilde{u}(t,x)=\zeta(t,x) v(t,x)$ for $t\geq S$. Then
$u:=\tilde{u}+(1-\zeta)v$ is obviously a solution we are after.
Its uniqueness follows immediately from Theorem 
\ref{theorem 2.4.2}. The theorem is proved.

The above proof allows one to get   corresponding
 estimates for the solution.
 We leave this to the interested reader.

{\bf Proof of Theorem \ref{theorem 2.7.1}}. 
Take
$\gamma_{b}$, $\gamma_{a}$, and $\lambda_{0}$
from Theorem \ref{theorem 2.4.2}, $f\in \cL_{p}(\bR^{d})$,  
$\lambda\geq\lambda_{0}$, and 
consider the equation
 \begin{equation}
                                             \label{2.15.3}
 \partial_{t} v+Lv-\lambda v=e^{- t}f
\end{equation}
in $\bR^{d+1}_{0}$.  As we have pointed out above,
we have $q_{1}< q$. Therefore, Assumption
\ref{assumption 2.7.2} (iii) implies that
Assumption \ref{assumption 2.9.1} is satisfied
for equation \eqref{2.15.3}. Other assumptions stated before 
Theorem \ref{theorem 2.4.2} are obviously satisfied
too. Hence, by Theorem 
\ref{theorem 2.9.1} equation \eqref{2.15.3} admits
a unique solution $v\in\cW^{2}_{p}(0)$. One easily checks that
for any $s\geq0$ the function $v(s+t,x)e^{s}$
as a function of $(t,x)\in\bR^{d+1}_{0}$
 also satisfies \eqref{2.15.3}. By
uniqueness $v(s+t,x)e^{s}=v(t,x)$, which implies that
$v(s,x)=e^{-s}u(x)$, where $u\in W^{2}_{p}(\bR^{d})$.
After that \eqref{2.15.3} is written as $Lu-(\lambda+1) u=f$.
This proves the existence in Theorem \ref{theorem 2.7.1}.
To prove uniqueness and estimate \eqref{2.4.01}
(with $\lambda+1$ in place of $\lambda$)
it suffices to introduce $v(t,x)=e^{-t}u(x)$,  observe that 
$v$ satisfies \eqref{2.15.3}, and use Theorem \ref{theorem 2.4.2}.
The theorem is proved.

\mysection{Auxiliary results}
                                               \label{section 2.15.3}

To emphasize which $b$ is used in the definition
of the operator $L$, write $L=L_{b}$.

\begin{lemma}
                                             \label{lemma 2.4.1}

There exist   constants  
$$
\gamma_{a}=\gamma(d,\delta,K,p)>0,\quad
 N=N(d,\delta,K,p,\varepsilon_{0}),
$$
$$
 \lambda_{0}=\lambda_{0}(d,\delta,K,p,\varepsilon_{0})\geq0
$$ 
such that, if Assumption \ref{assumption 2.7.1} (i) and
Assumption
\ref{assumption 2.4.1} ($\gamma_{a}$)  
are satisfied, then  
for any $T\in[-\infty,\infty)$,  $u \in 
\cW^{ 2}_{p}(T)$,
  $\lambda\geq\lambda_{0}$,  and 
 any $\bR^{d}$-valued
locally integrable 
to the power $p/(p-1)$  function
$\bar{b}=\bar{b}(t)$ on $\bR$
  we have
\begin{equation}
                                             \label{2.7.5}
\lambda\|u\|_{\cL_{p}(\bR^{d+1}_{T})}+
\|D^{2}u\|_{\cL_{p}(\bR^{d+1}_{T})}
\leq N\|L_{\bar{b}} u+\partial_{t}u
-\lambda u\|_{\cL_{p}(\bR^{d+1}_{T})}.
\end{equation}
 
\end{lemma}

Proof. First assume that $\bar{b}\equiv0$.
Since the coefficients of $L_{0}$ are bounded 
and $u \in 
\cW^{ 2}_{p}(T)$, the right-hand
side of  \eqref{2.7.5} is infinite unless $\partial_{t}u
\in\cL_{p}(\bR^{d+1}_{T})$, that is unless
$u\in W^{1, 2}_{p}(\bR^{d+1}_{T})$. In that case
our assertion is true 
by Theorem 6.4.1 and Remark 6.3.1 of \cite{book}.

In the case of general $\bar{b}$ take $u \in 
\cW^{ 2}_{p}(T)$ and introduce
$$
B(t)=\int_{0}^{t}\bar{b}(s)\,ds,\quad v(t,x)=u(t,x+B(t)),\quad
f=L_{\bar{b}} u+\partial_{t}u-\lambda u.
$$  
As is easy to see, the function $|\bar{b}(t)|\,|Du(t,x+B(t))|$
is locally summable in $\bR^{d+1}_{T}$ so that
 $v \in \cW^{ 2}_{p}(T)$ and
$$
\partial_{t}v(t,x)+[a^{ij}(t,x+B(t))D_{ij} 
-(\lambda+c(t,x+B(t))]v(t,x)
$$
$$
=f(t,x+B(t))=:g(t,x).
$$
By the above
$$
\lambda\|v\|_{\cL_{p}(\bR^{d+1}_{T})}+
\|D^{2}v\|_{\cL_{p}(\bR^{d+1}_{T})}
\leq N\|g\|_{\cL_{p}(\bR^{d+1}_{T})},
$$
which immediately yields \eqref{2.7.5}.
The lemma is proved.
\begin{remark}
In \cite{book} the assumption corresponding
to Assumption \ref{assumption 2.4.1} is much
weaker since in the corresponding
counterpart of \eqref{2.9.5}
there is no supremum over $x\in\bR^{d}$.
We need  our stronger assumption because
we need $a^{ij}(t,x+B(t))$ to satisfy
the assumption in \cite{book} for any function $\bar{b}$.

\end{remark}
\begin{remark}
                                       \label{remark 2.11.1}
The above proof
and the results in \cite{book} also show that for any $f\in\cL_{p}(\bR^{d+1}_{T})$
there exists a solution $u\in\cW^{2}_{p}(T)$ of the equation
$L_{\bar{b}} u+\partial_{t}u
-\lambda u=f$. Since the solution has the form $v(t,x-B(t))$
with $v\in W^{1,2}_{p}(\bR^{d+1}_{T})$, generally,
$\partial_{t}u$ is only locally summable in $t$.
\end{remark}

\begin{lemma}
                                        \label{lemma 2.11.1}
Suppose that Assumptions \ref{assumption 2.7.1} and
\ref{assumption 2.8.1} ($\gamma_{b}$) are satisfied
and let $n\in\{1,2,...\}$.
Then one can find a nonnegative function
$\xi\in C^{\infty}_{0}(B_{\alpha})$ which integrates to one
and a constant
$\beta_{n}=\beta(n,\gamma_{b},d,\alpha)$
such that, for almost any $t$, we have $|D^{n}\bar{b}_{\alpha}(t,x)|\leq
\beta_{n}$ on $\bR^{d }$, where $D^{n}\bar{b}$ is any derivative
of $\bar{b}_{\alpha}$ of order $n$ with respect to $x$ and
$$
\bar{b}_{\alpha}(t,x )= \int_{B_{\alpha}}
b(t,x-y)\xi(y)\,dy=\int_{\bR^{d}}b(t,y)\xi(x-y)\,dy.
$$

\end{lemma}

Proof. Take and fix any
nonnegative function
$\eta\in C^{\infty}_{0}(B_{1})$ which integrates to one
and set $\xi(x)=\alpha^{-d}\eta(x/\alpha)$. Due to Assumption
\ref{assumption 2.7.1} for almost any $t$ the function
$b(t,\cdot)$ is locally integrable on $\bR^{d}$ and hence
(for almost any $t$) the function
$\bar{b}_{\alpha}$ is well defined and infinitely differentiable
with respect to $x$.
Observe that
$$
D^{n}\bar{b}_{\alpha}(t,x)=\int_{B_{\alpha}}
b(t,x-y)D^{n}\xi(y)\,dy
=\int_{B_{\alpha}}
(b(t,x-y)-\bar{b}_{\alpha}(t,x))D^{n}\xi(y)\,dy.
$$
It follows that
$$
|D^{n}\bar{b}_{\alpha}(t,x)|\leq N
\int_{B_{\alpha}}
|b(t,x-y)-\bar{b}_{\alpha}(t,x)|\,dy
$$
$$
=N
\int_{B_{\alpha}}
\big|b(t,x-y)-\int_{B_{\alpha}}
b(t,x-z)\xi(z)\,dz\big|\,dy
$$
$$
=N
\int_{B_{\alpha}}
\big|\int_{B_{\alpha}}[b(t,x-y)-
b(t,x-z)]\xi(z)\,dz \,dy
$$
$$
\leq
N
\int_{B_{\alpha}}
 \int_{B_{\alpha}}|b(t,x-y)-
b(t,x-z)|\,dz \,dy
$$
and to get our assertion it only remains
to use H\"older's inequality. The lemma is proved.

\begin{corollary}
                                          \label{corollary 2.11.1}
Under Assumptions \ref{assumption 2.7.1} and 
\ref{assumption 2.8.1} ($\gamma_{b}$)
there exists a locally integrable 
 to the power $p/(p-1)$ function $K(t)$
on $\bR$ such that, for almost any $t$, we have on $\bR^{d}$ that 
\begin{equation}
                                             \label{2.8.2}
g(t,x):=|\bar{b}_{\alpha}(t,x)|\leq K(t)(1+|x|)
\end{equation}
  
\end{corollary}

Indeed by Lemma \ref{lemma 2.11.1} we have
$g(t,x) \leq g(t,0)   +\beta_{1}|x|$ and 
from Assumption \ref{assumption 2.7.1} (ii) we know that
$g(t,0)$ is locally integrable 
 to the power $p/(p-1)$.

\mysection{Proof of Theorem  \protect\ref{theorem 2.4.2}}
                                               \label{section 2.15.4}

We split the proof into several steps.

{\em Step 1. Introducing cut-off functions  
with time-dependent support}.
  First we take some $\gamma_{b}>0$  to be specified later,
suppose that Assumption \ref{assumption 2.8.1} ($\gamma_{b}$)
is satisfied with some $\alpha>0$,
 take $\beta=\beta_{1}$ from Lemma
\ref{lemma 2.11.1}, assume without loss
of generality that $\beta\geq1$,
and take $u\in\cW^{2}_{p}(0)$ such that
$u(t,x)=0$ for $t\geq\beta^{-1}$. Next,
fix  a  nonnegative
$ \zeta\in C^{\infty}_{0}(\bR^{d})$ with support in $B_{\alpha}$
and such that
\begin{equation}
                                             \label{2.4.3}
\int_{B_{\alpha}} \zeta^{p}(x) \,dx=1.
\end{equation}
 Also take  a point
$ x_{0} \in$  $  \bR^{d } $ and introduce
$x(t)=x_{x_{0}}(t)$ as a solution   of the problem
\begin{equation}
                                             \label{2.8.1}
x(t)=x_{0}+\int_{0}^{t}\bar{b}_{\alpha}
(s,x(s))\,ds,\quad t \in\bR,
\end{equation}
where $\bar{b}_{\alpha}$ is introduced in Lemma
\ref{lemma 2.11.1}.
Owing to Lemma
\ref{lemma 2.11.1} and Corollary \ref{corollary 2.11.1}
equation  \eqref{2.8.1} admits a  unique 
solution which is infinitely differentiable
with respect to $x_{0}$ because $\bar{b}_{\alpha}(t,x)$
is infinitely differentiable in $x$.

Set
$$
  \bar{b}_{x_{0}}(t)=\bar{b}_{\alpha}(t,x_{x_{0}}(t)).
$$
$$
L_{x_{0}} =a ^{ij}(t,x)D_{ij}+
\bar{b}_{x_{0}}(t)D_{i}-c (t,x), 
$$
$$
\eta_{x_{0}}(t,x)=\zeta(x-x_{x_{0}}(t)), 
\quad v_{x_{0}}(t,x)= u(t,x)  \eta_{x_{0}}(t,x),
$$
$$
f:=L u+\partial_{t}u -\lambda u.
$$

Observe that
$$
\partial_{t}\eta_{x_{0}}(t,x)+\bar{b}^{i}_{x_{0}}(t)
D_{i} \eta_{x_{0}} (t,x)=0, 
$$
which implies that
$$
\partial_{t}v_{x_{0}}+L_{x_{0}}v_{x_{0}}
-\lambda v_{x_{0}} =
\eta_{x_{0}}(\partial_{t}u+L_{x_{0}}u
-\lambda u)
$$
$$+
u (\partial_{t}\eta_{x_{0}}+L_{x_{0}}
\eta_{x_{0}}+c \eta_{x_{0}})
+2a^{ij}(D_{i}\eta_{x_{0}})D_{j}u 
$$
\begin{equation}
                                             \label{2.4.4}
=\eta_{x_{0}}  f -f^{1}_{x_{0}}
+f^{2}_{x_{0}}+f^{3}_{x_{0}},
\end{equation}
where
$$
f^{1}_{x_{0}}=\eta_{x_{0}} 
(b^{i}-\bar{b}^{i}_{x_{0}})
 D_{i}u,\quad f^{2}_{x_{0}}=
 u a^{ij} D_{ij}\eta _{x_{0}},
$$
$$
f^{3}_{x_{0}}= 2a^{ij}(D_{i}\eta_{x_{0}})D_{j}u.
$$

{\em Step 2.  Estimating the right-hand side of \eqref{2.4.4}}.
Observe that if   $\eta_{x_{0}}(t,x)\ne0$, then
$|x-x_{x_{0}}(t)|\leq\alpha$ and we may certainly assume that
$\zeta^{p}\leq N(d)\alpha^{-d}$, so that
$$
\|f^{1}_{x_{0}}\|_{\cL_{p}(\bR^{d+1}_{0})}^{p}
=\int_{0}^{\infty}\int_{B_{\alpha}+x_{x_{0}}(t) }
|(b^{i}-\bar{b}^{i}_{x_{0}})
 \eta_{x_{0}}D_{i}u|^{p}\,dxdt
$$
$$
\leq N(d)\alpha^{-d}\int_{0}^{\infty}I(t)\,dt,
$$
where
$$
I(t)=\int_{B_{\alpha}+x_{x_{0}}(t) }
|(b^{i}-\bar{b}^{i}_{x_{0}})
  D_{i}u|^{p}(t,x)\,dx .
$$

If $p\leq d$ we use
 H\"older's 
inequality,
and embedding theorems to obtain
$$
I(t)\leq \big(\int_{B_{\alpha}+x_{x_{0}}(t) }
|b -\bar{b} _{x_{0}}|^{q}\,dx\big)^{\frac{p}{ q}}
\big(\int_{B_{\alpha}+x_{x_{0}}(t) }
|D  u |^{\frac{pq}{  q-p}}\,dx\big)^{
\frac{ q-p}{ q}} 
$$
$$
\leq N\mu^{\frac{p}{ q}}(t,x_{x_{0}}(t),x_{0})
[v(t,x_{x_{0}}(t))+w(t,x_{x_{0}}(t))],
$$
where $N=N(d,p )$, 
$$
 \mu(t,y,x_{0}):=\int_{B_{\alpha}+y }
|b(t,x) -\bar{b} _{x_{0}}(t)|^{q}\,dx,
$$
$$
v(t,y)=\alpha^{p-p'}\int_{B_{\alpha}+y }
|D^{2}u(t,x)|^{p}\,dx,\quad w(t,y)
=\alpha^{-p-p'}\int_{B_{\alpha}+y }
| u(t,x)|^{p}\,dx,
$$
and $p'=pd/q$ ($\leq p$).

We also note that $\xi\leq N(d)\alpha^{-d}$ (see
Lemma \ref{lemma 2.11.1})
and  by Assumption \ref{assumption 2.3.1} ($\gamma_{b}$)
we find that
$$
 \mu(t,y,x_{0})=\int_{B_{\alpha}+x_{x_{0}}(t) }
\big|\int_{B_{\alpha}+x_{x_{0}}(t) }[
b(t,x) -b(t,y)]\xi( x_{x_{0}}(t)-y)\,dy\big|^{q}\,dx
$$
$$
\leq N(d)\alpha^{-d}\int_{B_{\alpha}+x_{x_{0}}(t) }
 \int_{B_{\alpha}+x_{x_{0}}(t) }|
b(t,x) -b(t,y) |^{q}\,dxdy\leq N(d )\gamma_{b}.
$$
Hence
$$
I(t) \leq N\gamma_{b}^{p/q}
[v(t,x_{x_{0}}(t))+w(t,x_{x_{0}}(t))],
$$
where $N=N(d,p )$.  
 This
estimate   also holds if $p>d$,
which is seen if we start like
$$
I(t)\leq
\int_{B_{\alpha}+x_{x_{0}}(t) }
|b -\bar{b} _{x_{0}}|^{p}\,dx\sup_{
B_{\alpha}+x_{x_{0}}(t) }
|D  u |^{p}. 
$$
Thus,
$$
\|f^{1}_{x_{0}}\|_{\cL_{p}(\bR^{d+1}_{0})}^{p}
\leq N(d,p)\gamma_{b}^{p/q}\alpha^{-d}\int_{0}^{\infty}
[v(t,x_{x_{0}}(t))+w(t,x_{x_{0}}(t))]\,dt.
$$

The following   estimates of $f^{2}$ and $f^{3}$
are straightforward:
$$
\| f^{2}_{x_{0}}\|_{\cL_{p}(\bR^{d+1}_{0})}^{p}
\leq N \int_{0}^{\infty}\int_{\bR^{d}
}I_{B_{\alpha}}( x_{x_{0}}(t)-x) |u(t,x)|^{p}\,dxdt,
$$
$$
\| f^{3}_{x_{0}}\|_{\cL_{p}(\bR^{d+1}_{0})}^{p}
\leq N\int_{0}^{\infty}\int_{\bR^{d}
}I_{B_{\alpha}}( x_{x_{0}}(t)-x)|Du(t,x)|^{p}\,dxdt,
$$
where $N=N(d,p,K,\alpha)$.
 
Now, provided that $\lambda\geq
\lambda_{0}$ with $\lambda_{0}$ from Lemma \ref{lemma 2.4.1},
 equation \eqref{2.4.4} and
  Lemma \ref{lemma 2.4.1} lead to  
$$
\lambda^{p}\|u\eta_{x_{0}}\|_{\cL_{p}(\bR^{d+1}_{0})}^{p}+\|
D^{2}(u\eta_{x_{0}})\|_{\cL_{p}(\bR^{d+1}_{0})}^{p}
\leq N_{1} \|f\eta_{x_{0}}\|_{\cL_{p}(\bR^{d+1}_{0})}^{p}
$$
$$
+N_{1}\gamma_{b}^{p/q}\alpha^{-d}\int_{0}^{\infty}
[v(t,x_{x_{0}}(t))+w(t,x_{x_{0}}(t))]\,dt
$$
$$
+N_{2}\int_{0}^{\infty}\int_{\bR^{d}
}I_{B_{\alpha}}( x_{x_{0}}(t)-x)( |u(t,x)|^{p}
+|Du(t,x)|^{p}) \,dxdt,
$$
where and below by $N_{1}$ we denote
generic constants depending only on $ d,\delta,K,p,\varepsilon_{0} $
and by $N_{2}$ constants depending
 only on the same parameters and $\alpha$.
By writing what $D^{2}(u\eta_{x_{0}})$ is, we conclude
 $$
\lambda^{p}\|u\eta_{x_{0}}\|_{\cL_{p}(\bR^{d+1}_{0})}^{p}+\|
\eta_{x_{0}}D^{2}u\|_{\cL_{p}(\bR^{d+1}_{0})}^{p}
\leq N_{1} \|f\eta_{x_{0}}\|_{\cL_{p}(\bR^{d+1}_{0})}^{p}
$$
$$
+N_{1}\gamma_{b}^{p/q}\alpha^{-d}\int_{0}^{\infty}
[v(t,x_{x_{0}}(t))+w(t,x_{x_{0}}(t))]\,dt
$$
\begin{equation}
                                              \label{2.5.1}
+N_{2}\int_{0}^{\infty}\int_{\bR^{d}
}I_{B_{\alpha}}( x_{x_{0}}(t)-x)( |u(t,x)|^{p}
+|Du(t,x)|^{p}) \,dxdt.
\end{equation}

{\em Step 3.  Integrating through \eqref{2.5.1}  with respect to
$x_{0}$}. One knows that for each   $t$,
the mapping $x_{0}\to x_{x_{0}}(t)$ is a diffeomorphism
with Jacobian determinant given by
$$
\bigg|\frac{\partial x_{x_{0}}(t)}{
\partial x_{0}}\bigg|(x_{0})=\exp\int_{0}^{t}\sum_{i=1}^{d}[D_{i}
\bar{b}_{\alpha}^{i}](s,x_{x_{0}}(s))\,ds.
$$
By Lemma \ref{lemma 2.11.1}
$$
e^{-N\beta t}\leq \bigg|\frac{\partial x_{x_{0}}(t)}{
\partial x_{0}}\bigg|(x_{0}) \leq e^{N\beta t},
$$
where   $N$ depends only on $d$.
  Therefore, for any nonnegative
Lebesgue measurable function $w(x)$ we have
$$
e^{-N\beta t}
\int_{\bR^{d}}w(y)\,dy\leq
\int_{\bR^{d}}w(x_{x_{0}}(t))\,dx_{0}\leq e^{N\beta t}
\int_{\bR^{d}}w(y)\,dy .
$$
In particular, since
$$
\int_{\bR^{d}}|\eta_{x_{0}}(t,x)|^{p}\,dx_{0}=
\int_{\bR^{d}}|\zeta(x-x_{x_{0}}(t )|^{p}\,dx_{0} ,
$$
we have
$$
e^{-N\beta t}=e^{-N\beta t}
 \int_{\bR^{d}}|\zeta(x-y)|^{p}\,dy 
 \leq\int_{\bR^{d}}|\eta_{x_{0}}(t,x)|^{p}\,dx_{0}
\leq e^{N\beta t} ,
$$
so that
$$
\int_{0}^{\infty}\int_{\bR^{d}}
[v(t,x_{x_{0}}(t))+w(t,x_{x_{0}}(t))]\,dx_{0}dt
$$
$$
\leq \int_{0}^{\infty}e^{ N\beta t}\int_{\bR^{d}}
[v(t,y)+w(t,y)]\,dydt
$$
$$
=N(d)\alpha^{d}\int_{0}^{\infty}e^{ N\beta t}\int_{\bR^{d}}
[\alpha^{p-p'}|D^{2}u|^{p}+\alpha^{-p-p'}|u|^{p}
(t,x)\,dxdt.
$$
Similarly one treats other terms in \eqref{2.5.1}.
For instance,
$$
\int_{\bR^{d}} 
 \|f\eta_{x_{0}}\|_{\cL_{p}(\bR^{d+1}_{0})}^{p}\,dx_{0}
\leq\int_{0}^{\infty}e^{ N\beta t}
\int_{\bR^{d}}|f(t,x)|^{p}\,dxdt,
$$
$$
\int_{\bR^{d}} 
 \|u\eta_{x_{0}}\|_{\cL_{p}(\bR^{d+1}_{0})}^{p}\,dx_{0}
\geq\int_{0}^{\infty}e^{ -N\beta t}
\int_{\bR^{d}}|u(t,x)|^{p}\,dxdt.
$$
We also observe that we need not integrate with respect to  $t$
beyond $\beta^{-1}$ which allows us to conclude
from  \eqref{2.5.1} that 
 $$
\lambda^{p}\|u \|_{\cL_{p}(\bR^{d+1}_{0})}^{p}+\|
 D^{2}u\|_{\cL_{p}(\bR^{d+1}_{0})}^{p}
\leq N_{1} \|f \|_{\cL_{p}(\bR^{d+1}_{0})}^{p}
$$
$$
+N_{1}\gamma^{p/q}_{b}\alpha^{p-p'}\|
D^{2}u\|_{\cL_{p}(\bR^{d+1}_{0})}^{p}
+N_{2}  \|D
 u\|_{\cL_{p}(\bR^{d+1}_{0})}^{p}
+N_{2}(1+\gamma^{p/q}_{b}) \|
 u\|_{\cL_{p}(\bR^{d+1}_{0})}^{p},
$$
where still $N_{1}=N_{1}(d,\delta,K,p,\varepsilon_{0})$
and $N_{2}$ depends only on the same parameters and $\alpha$.

We now can specify $\gamma_{b}(d,\delta,K,p,\varepsilon_{0})$. 
We take it in such a way
that $N_{1}\gamma^{p/q}_{b}\leq1/2$. Then
(recall that $\alpha\in(0,1]$ and $p\geq p'$) we obtain
$$
\lambda^{p}\|u \|_{\cL_{p}(\bR^{d+1}_{0})}^{p}+\|
 D^{2}u\|_{\cL_{p}(\bR^{d+1}_{0})}^{p}
\leq N_{1} \|f \|_{\cL_{p}(\bR^{d+1}_{0})}^{p}
$$
\begin{equation}
                                              \label{2.5.01}
+N_{2}  \|D
 u\|_{\cL_{p}(\bR^{d+1}_{0})}^{p}
+N_{2}  \|
 u\|_{\cL_{p}(\bR^{d+1}_{0})}^{p}.
 \end{equation}
Now we show  how to choose $\lambda_{0}
(d,\delta,K,p,\varepsilon_{0},\alpha)$. We take it
larger than the one from Lemma
\ref{lemma 2.4.1} and such that for $\lambda\geq
\lambda_{0}$ interpolation inequalities would allow
us to absorb the last two terms in \eqref{2.5.01}
into its left-hand side. Then for $\lambda\geq
\lambda_{0}$ we get
 $$
\lambda^{p}\|u \|_{\cL_{p}(\bR^{d+1}_{0})}^{p}+\|
 D^{2} u \|_{\cL_{p}(\bR^{d+1}_{0})}^{p}
\leq N \|f \|_{\cL_{p}(\bR^{d+1}_{0})}^{p}
$$
with $N=N(d,\delta,K,p,\varepsilon_{0})$, 
provided that $u\in\cW^{2}_{p}(0)$ is such that
$u(t,x)=0$ for $t\geq\beta^{-1}$.

{\em Step 4.  Case of $u$ not compactly supported in $t$}.
To pass to the general case we take a nonnegative
 function $\chi\in C^{\infty}_{0}(0,1)$ such that
$$
\int_{0}^{1}\chi^{p}(t)\,dt=1
$$
and set $\kappa_{t_{0}}(t)=\beta^{1/p}\chi(\beta t-
\beta t_{0})$. Recall that $\gamma_{b}$ is fixed above,
so that in Lemma \ref{lemma 2.11.1}
we have $\beta_{1}=\beta=\beta(d,\delta,K,p,\varepsilon_{0},\alpha)$.
Next, for each $t_{0}$
we have that $u_{t_{0}}(t,x):=u(t,x)\kappa_{t_{0}}(t)$
belongs to $\cW^{2}_{p}(T\vee t_{0})$ and $u_{t_{0}}(t,x)=0$
for $t\geq T\vee t_{0}+\beta^{-1}$. The result of
the above particular case implies that
$$
\int_{T}^{\infty}\int_{\bR^{d}}
[\lambda^{p}|u|^{p}+|D^{2}u|^{p}](t,x)
\kappa^{p}_{t_{0}}(t)\,dxdt
$$
$$
\leq N_{1}\int_{T}^{\infty}\int_{\bR^{d}}
[|f|^{p}\kappa^{p}_{t_{0}}+| u|^{p}(\kappa'_{t_{0}})^{p}](t,x)\,dxdt,
$$
provided that $\lambda\geq\lambda_{0}(d,\delta,K,p,
\varepsilon_{0},\alpha)$.
By integrating through this with respect to $t_{0}$
over $\bR$ we obtain
 $$
\lambda^{p}\|u \|_{\cL_{p}(\bR^{d+1}_{T})}^{p}+\|
 D^{2} u \|_{\cL_{p}(\bR^{d+1}_{T})}^{p}
\leq N _{1}\|f \|_{\cL_{p}(\bR^{d+1}_{T})}^{p}
+N_{2}\|u \|_{\cL_{p}(\bR^{d+1}_{T})}^{p} .
$$
  Now it only remains to increase $\lambda_{0}$ if needed
to absorb the last term on the right into the left-hand side.

The theorem is proved.

\mysection{Proof of   Theorem~\protect\ref{theorem 2.9.1}} 
                                               \label{section 2.15.5}
We will only concentrate on the existence
since uniqueness follows from \eqref{2.4.1}.
We start with the following.
\begin{lemma}
                                     \label{lemma 2.10.1}
Let $b$ be a measurable function
and let Assumption \ref{assumption 2.9.1} be satisfied.
Let $b_{n}(t,x)$, $n=1,2,...$, be $\bR^{d}$-valued 
measurable functions
on $\bR^{d+1}_{T}$ such that  
\begin{equation}
                                                \label{2.10.5}
\int_{s}^{t}\big(\int_{B_{R}}|b_{n}(\tau,x)-
b (\tau,x)|^{q_{1}}
\,dx\big)^{r_{1}}\,d\tau\to0
\end{equation}
as long as $s,t\in\bR$, $s<t$, and $R\in(0,\infty)$.
Finally, let $u$, $u_{n}$, $Du$, $Du_{n}$,
$D^{2}u$, $D^{2}u_{n}\in\cL_{p}
(\bR^{d+1}_{T})$, $n=1,2,...$, and assume that
\begin{equation}
                                                \label{2.10.10}
(u_{n},Du_{n},D^{2}u_{n})\to(u ,Du ,D^{2}u )
\end{equation}
weakly in $\cL_{p}
(\bR^{d+1}_{T})$. 

Then $b^{i}_{n}D_{i}u_{n}
\to b^{i} D_{i}u $
  in the sense of distributions on $\bR^{d+1}_{T}$
and $b^{i} D_{i}u$ is locally summable on $\bR^{d+1}_{T}$.
\end{lemma}

Proof. Take $\kappa>0$, $T<s<t$, $R>0$, and 
$\phi,\psi\in C^{\infty}_{0}(\bR^{d+1}_{T})$ with support
in $Q:=(s,t)\times B_{R}$  such that
\begin{equation}
                                                \label{2.10.03}
\int_{s}^{t}\big(\int_{B_{R}}|b -\psi|^{q_{1}}
(\tau,x)\,dx\big)^{r_{1}}\,d\tau\leq\kappa.
\end{equation}
Use the notation $(g,h)$ for the integral of $gh$
over $\bR^{d+1}_{T}$ and write
$$
I_{n }:=|(\phi,b^{i}_{n}D_{i}u_{n})
- (\phi,b^{i} D_{i}u )|
\leq(|\phi|,|b _{n}-b  |
\,|D u_{n}|)
$$
$$
+(|\phi|,|b  -\psi|
\,|D u_{n}|)+(|\phi|,|b  -\psi|
\,|D u |)+(\phi\psi^{i},D_{i}
u_{n}-D_{i}u ).
$$
Here the last term goes to zero since $D_{i}
u_{n}\to D_{i}u $ weakly.
To estimate the remaining terms on the right we use
embedding theorems. For instance,
$$
(|\phi|,|b -\psi|
\,|D u_{n}|)\leq N\int_{s}^{t}\int_{B_{R}}
|b -\psi|(\tau,x)|Du_{n}(\tau,x)|
\,dxd\tau
$$
$$
\leq N\int_{s}^{t}J_{1}(\tau)J_{2,n}(\tau)\,d\tau
\leq\big(\int_{s}^{t}J_{1}^{p/(p-1)}(\tau)\,d\tau\big)
^{(p-1)/p}
\big(\int_{s}^{t}J_{2,n}^{p}(\tau)\,d\tau)^{1/p},
$$
where
$$
J_{1}(\tau)=\big(\int_{B_{R}}
|b -\psi|^{q_{1}}(\tau,x)\,dx\big)^{1/q_{1}}
$$
 $$
J_{2,n}(\tau)=\big(
\int_{B_{R}}|Du_{n}(\tau,x)|^{q_{1}/(q_{1}-1)}
\,dx\big)^{(q_{1}-1)/q_{1}}.
$$

In light of \eqref{2.10.03} and the fact that $p/[q_{1}(p-1)]
=r_{1}$ we have
$$
\int_{s}^{t}J_{1}^{p/(p-1)}(\tau)\,d\tau\leq\kappa.
$$
Furthermore, by embedding theorems and the fact that
$q_{1}/(q_{1}-1)=pd/(d-p)$ for $p<d$ and $q_{1}/(q_{1}-1)<\infty$
in any case, we obtain
$$
J_{2,n}(\tau)\leq N(\|u_{ n}(\tau,\cdot)\|_{\cL_{p}(\bR^{d})}
+\|D^{2}u_{ n}(\tau,\cdot)\|_{\cL_{p}(\bR^{d})}),
$$
so that
$$
\nlimsup_{n\to\infty}\int_{s}^{t}J_{2,n}^{p}(\tau)\,d\tau\leq
N\nlimsup_{n\to\infty}(\|u_{ n}\|_{\cL_{p}(\bR^{d+1}_{T })}^{p}
+ \|D^{2}u_{ n}\|_{\cL_{p}(\bR^{d+1}_{T })}^{p})
$$
which  is finite due to the assumed weak convergence of $u_{n}$
and $D^{2}u_{n}$.

Similarly one estimates $(|\phi|,|b -\psi|
\,|D u |)$ and $(|\phi|,|b _{n}-b  |
\,|D u_{n}|)$ invoking \eqref{2.10.5}
 in the case of the latter. Hence 
$$
\nlimsup_{n }I_{n }\leq N\kappa^{(p-1)/p}
$$
with $N$ independent of $\kappa$
and, since $\kappa>0$ is arbitrary, $I_{n }\to0$ as $n \to
\infty$. The arbitrariness of $\phi$ finishes proving our
claim.

The reader might have noticed that the above computations
are only valid if we knew
that $\phi b^{i}D_{i}u$ and $\phi b^{i}_{n}D_{i}u_{n}$
are  summable   at least for large $n$. To close this gap
it suffices to fix $n$ and do
the above estimates with $b_{n}=0$, $u_{n}=0$.
This would prove that $\phi b^{i}D_{i}u$
is  summable. Due to  \eqref{2.10.5} the functions
$b_{n}$ satisfy \eqref{2.10.1} for all large $n$
and, as for $\phi b^{i}D_{i}u$, this implies that
$\phi b^{i}_{n}D_{i}u_{n}$
are  summable    for large $n$.
The lemma is proved.

{\bf  Proof of Theorem \ref{theorem 2.9.1}}.
Recall that $\gamma_{a},\gamma_{b}$, and $\lambda_{0}$
are taken
from Theorem \ref{theorem 2.4.2}, for $n=1,2,...$ define
$\kappa_{n}(t)=(-n)\vee t\wedge n$ and $b^{i}_{n}
= \kappa_{n}(b^{i})$, and observe that, since the $b_{n}$
are bounded for each $n$, by Theorem 6.4.1 of \cite{book},
for any $\lambda\geq\lambda_{0}$ and $f\in\cL_{p}(\bR^{d+1}_{T})$,
there exist $u_{n}\in W^{1,2}_{p}(\bR^{d+1}_{T})$ satisfying
\begin{equation}
                                              \label{2.10.9}
L_{b_{n}}u_{n}+\partial_{t}u_{n} -\lambda u_{n}=f
\end{equation}
in $\bR^{d+1}_{T}$.

Notice that $\kappa_{n}$ are Lipschitz continuous functions
with Lipschitz constant 1. It follows that $b_{n}$ satisfies
Assumption \ref{assumption 2.8.1} ($\gamma_{b}$)
(with the same $\gamma_{b}$). Hence by Theorem
\ref{theorem 2.4.2} 
\begin{equation}
                                                 \label{2.9.8}
\|u_{ n}\|_{\cL_{p}(\bR^{d+1}_{T })}
+\|D^{2}u_{ n}\|_{\cL_{p}(\bR^{d+1}_{T })}\leq
N<\infty,
\end{equation}
where $N$ is independent of   $n$.

Owing to \eqref{2.9.8} there is a subsequence $n'\to\infty$
and a function $u $ such that
$$
u_{ n'}\to u 
,\quad Du_{ n'}\to Du ,\quad
D^{2}u_{ n'}\to D^{2}u 
$$
weakly in $\cL_{p}(\bR^{d}_{T})$.
Consequently,
also weakly in $\cL_{p}(\bR^{d}_{T})$, we have
\begin{equation}
                                                \label{2.10.3}
L_{0}u_{ n'}\to L_{0}u .
\end{equation}

Next, $|b_{n}|\leq |b|$ and $b_{n}\to b$ as $n\to\infty$
(a.e.). By the dominated convergence theorem
we have that condition \eqref{2.10.5} is satisfied.
Then, according to Lemma \ref{lemma 2.10.1} we have
that  $b^{i}_{n}D_{i}u_{n}
\to b^{i} D_{i}u $
  in the sense of distributions on $\bR^{d+1}_{T}$
and $b^{i} D_{i}u$ is locally summable on $\bR^{d}_{T}$.
Now by passing to the limit in the sense of distributions
in \eqref{2.10.9} we see that $Lu+\partial_{t}u-\lambda u=f$
in $\bR^{d+1}_{T}$ and $\partial_{t}u$ is locally summable.
It follows that $u\in\cW^{2}_{p}(T)$ and the theorem is proved.

\mysection{An example}
                                               \label{section 2.15.6}
Let $d=1$ and consider the elliptic equation
\begin{equation}
                                           \label{2.12.1}
Lu:=u''-2bx u'- 2 u=-2f.
\end{equation}
We claim that if $b\geq p$, then there exist
functions $f\in\cL_{p}(\bR)$ such that equation
\eqref{2.12.1} does not have solutions $u\in W^{2}_{p}(\bR)$. This fact
is, actually, known from a very interesting article \cite{Me},
where the spectrum of the Ornstein-Uhlenbeck operator is found
in the multidimensional case in $\cL_{p}$ spaces. 
In the case of \eqref{2.12.1} the result of \cite{Me} says that
this equation is uniquely solvable in $W^{2}_{p}(\bR)$
iff $p>b$.  We give an independent and short proof
of our claim for completeness of presentation.
A very rough idea why this happens in our situation
 is that if the resolvent
operator of this equation were bounded in $\cL_{p}(\bR)$,
then its adjoint would also be bounded, but this
adjoint is the resolvent operator of the adjoint
equation that  has $c$ with a wrong sign.

We know (see, for instance,
\cite{KP}) that for any $f\in C^{\infty}_{0}(\bR)$
 there exists a unique
smooth and bounded solution $u$ of \eqref{2.12.1}.
It is well known that $u$ is given by
$$
u(x)=\int_{\bR^{d}}f(y)g(x,y)\,dy=:Rf(x),
$$
where
$$
g(x,y)=\int_{0}^{\infty}e^{-t}p(t,x,y)\,dt,
$$
$$
\sigma^{2}(t)=\int_{0}^{t}e^{-2bs}\,ds,\quad
p(t,x,y)= \frac{1}{\sqrt{2\pi\sigma^{2}(t)}}
\exp\big[-\frac{(y-xe^{-bt})^{2}}{2\sigma^{2}(t)}\big].
$$
By the maximum principle any solution of class $W^{2}_{p}(\bR)$
coincides with $Rf$. 
Now to prove our claim, it suffices to show that,
if $f\geq1 $ on $(0,1)$ and $f\geq0$ on $R$, then we have
$u=Rf\not\in\cL_{2}(\bR)$. Observe that
$$
u(x)\geq\int_{0}^{\infty}e^{-t}\big(\int_{0}^{1}
p(t,x,y)\,dy\big)\,dt
$$
and since for $t\geq1$ the function $\sigma^{2}(t)$
is bounded away from zero and infinity, we have
$$
u(x)\geq \varepsilon
\int_{1}^{\infty}e^{-t}\big(\int_{0}^{1}
\exp[-N(y-xe^{-bt})^{2}]\,dy\big)\,dt,
$$ 
where $\varepsilon>0$ and $N$ are some constants.
Obviously, the interior integral is bigger than a constant
$\varepsilon_{1}>0$ if $0\leq xe^{-bt}\leq1$, that is if
$x\geq0$ and
$t\geq b^{-1}\log x$.
Thus, for $x\geq e^{b}$
$$
u(x)\geq \varepsilon\varepsilon_{1}\int_{b^{-1}\log x}^{\infty}
e^{-t}\,dt=\frac{\varepsilon\varepsilon_{1}}{x^{1/b}},
$$
which is not in $\cL_{p}(e^{b},\infty)$ for $b\geq p$.


\begin{thebibliography}{mm}

\bibitem{BL} A. Bensoussan  and J.-L. Lions, ``Applications 
of variational inequalities in stochastic control", North-
Holland, Amsterdam, 1982.
 

\bibitem{CV} P. Cannarsa  and V. Vespri, {\em
 Existence and uniqueness results
for a nonlinear stochastic partial differential equation\/}, in 
Stochastic Partial Differential Equations and Applications
Proceedings, G. Da Prato and L. Tubaro (eds.), Lecture Notes in
Math., Vol. 1236, pp. 1-24, Springer Verlag, 1987.

\bibitem{CF} G. Cupini and S. Fornaro, {\em Maximal regularity in $L^
p(\bR^N)$ for a class of elliptic operators with unbounded
coefficients\/},  Differential Integral Equations, Vol.  17  (2004), 
 No. 3-4, 259-296. 

\bibitem{GL} M. Geissert and A. Lunardi,  {\em Invariant
measures and maximal $L\sp 2$ regularity for nonautonomous
Ornstein-Uhlenbeck equations\/},  J. Lond. Math. Soc. (2), Vol.
  77  (2008), 
No. 3, 719-740.

\bibitem{Gy93}   I. Gy\"ongy,
{\em Stochastic partial differential equations on
Manifolds, I\/}, Potential Analysis, Vol.  2 (1993), 101-113.

\bibitem{GK}   I. Gy\"ongy and N.V. Krylov,  {\em
  On stochastic partial differential equations with unbounded
 coefficients\/},  Potential Analysis, Vol. 1 (1992), No. 3, 233-256.

\bibitem{book} N.V. Krylov,
``Lectures on elliptic and parabolic equations
in Sobolev spaces", Amer.
Math. Soc., Providence, RI, 2008.

\bibitem{KP} N.V. Krylov and E. Priola, {\em
Elliptic and parabolic second-order PDEs with growing coefficients\/},
submitted to Comm in PDEs, http://arXiv.org/abs/0806.3100
  
\bibitem{LV} A. Lunardi and V. Vespri,
{\em Generation of strongly continuous semigroups by elliptic
operators with unbounded coefficients in $L^p(\bR^n)$\/},
Rend. Istit. Mat. Univ. Trieste 28 (1996), suppl., 
251-279 (1997).

\bibitem{Me} G. Metafune, {\em
$L^p$-spectrum of Ornstein-Uhlenbeck operators\/}, Annali della Scuola
Normale Superiore di Pisa - Classe di Scienze, S\'er. 4, Vol. 30,
 No. 1
(2001),  97-124.

\bibitem{MP} G.
Metafune, J. Pr\"uss,  R. Schnaubelt, and A. Rhandi,  {\em
$L\sp p$-regularity for elliptic operators with unbounded 
coefficients\/},
Adv. Differential Equations, Vol. 10 (2005), No. 10, 1131-1164. 

\bibitem{PR} J.
Pr\"uss, A. Rhandi,  and R. Schnaubelt, {\em
The domain of elliptic operators on $L^p(\bR^d)$ with
unbounded drift coefficients\/}, Houston J. Math., Vol. 32
(2006), No. 2, 563-576.
\end{thebibliography}
\end{document}